\documentclass[11pt,reqno]{amsart}
\usepackage{amsfonts,amssymb,diagrams}

\begin{document}
\title{A bordism approach to string topology}
\author[D. Chataur]
{David Chataur\\
CRM-Barcelona}

\thanks{The author was partially supported by DGESIC grant
PB97-0202.}

\date{\today}

\email{dchataur@crm.es}

\maketitle

\begin{abstract}
Using intersection theory in the context of Hilbert manifolds and
geometric homology we show how to recover the main operations of
string topology built by M. Chas and D. Sullivan. We also study
and build an action of the homology of reduced Sullivan's chord
diagrams on the singular homology of free loop spaces, extending
previous results of R. Cohen and V. Godin and unifying part of the
rich algebraic structure of string topology as an algebra over the
prop of these reduced diagrams. Some of these operations are
extended to spaces of maps from a sphere to a compact manifold.
\end{abstract}

\section{Introduction}

The study of spaces of maps is an important and difficult task of
algebraic topology. In this paper we study $n$-sphere spaces, they
are topological spaces of unbased maps from a $n$-sphere into a
manifold. Our aim is to study the algebraic structure of the
homology of these spaces. Let us begin by giving some motivations
for the study of such spaces. We focus on free loop spaces
($n=1$).
\\
\\
The study of free loop spaces over a compact oriented manifold
plays a central role in algebraic topology. There is a non
exhaustive list of topics where free loop spaces appear. Let us
review some of them:
\\
- Free loop spaces are one
of the main tool in order to study closed geodesics on Riemannian
manifolds. Let $M$ be Riemannian, compact, connected, simply connected,
of dimension greater than one. D. Gromoll and W. Meyer proved that
there exist infinitely many (geometrically distinct)
periodic geodesics on an arbitrary Riemannian manifold if the Betti numbers of the free loop space of
$M$ are unbounded \cite{GM}. And using methods of rational homotopy theory M. Vigu\'e and D. Sullivan
showed that these rational Betti numbers are unbounded if and only if the rational cohomology has at
least two generators \cite{VS}.
\\
- According to works of D. Burghelea and Z. Fiedorowicz \cite{BF}
of T. Goodwillie \cite{Go} and of J.D.S. Jones \cite{JJ}, the
cohomology of free loop spaces is strongly related to Hochschild
homology and the $S^1$-equivariant cohomology of free loop spaces
is also related to cyclic homology and over a field of
characteristic zero to Waldhausen algebraic $K$-theory. Hence
rationally, the homotopy of free loop spaces is related to the
space of automorphisms of manifolds (via Waldhausen's theory). To
go further in this direction there are also some relations between
the suspension spectrum of the free loop spaces and topological
cyclic homology \cite{BHM}.
\\
- Analysis on such spaces has always been a source of constant
inspiration, for example one can cite the work of E. Witten on
Index of Dirac operators on Free loop spaces \cite{W}, this was
the beginning of the theory of elliptic genera and elliptic
cohomology \cite{S}. This shed new lights on the periodicity
phenomena in stable homotopy and the analysis underlying it
\cite{AM}.
\\
- One can also cite the work of K. T. Chen who has introduced a
chain complex based on iterated integrals that computes the
cohomology of free loop spaces over a manifold \cite{Ch}. Iterated
integrals give a De-Rham theory for path spaces. In particular,
this theory has found some applications in the algebraic
interpretation of index theory on free loop spaces by means of
cyclic homology \cite{GJP1} and \cite{GJP2}. Let us notice that
the theory of iterated integrals s also related to algebraic
geometry (see \cite{H} for a survey).
\\
\\
More recently the discovery by M. Chas and D. Sullivan of a
Batalin Vilkovisky structure on the singular homology of these
spaces \cite{CS} had a deep impact on the subject and has revealed
a part of a very rich algebraic structure \cite{CS2}, \cite{CG}.
The $BV$-structure consists in
\\
- A loop product $-\bullet -$ which is commutative and
associative, it can be understood as an attempt to perform
intersection theory of families of closed curves,
\\
- A loop bracket $\{-,-\}$, half of this bracket controls the
commutativity up to homotopy at the chain level of the loop
product,
\\
- An operator $\Delta$ coming from the action of $S^1$ on the free loop space ($S^1$ acts by
reparametrization of the loops).
\\
M. Chas and D. Sullivan use (in \cite{CS}) "classical intersection
theory of chains in a manifold" . This structure has also been
defined in a purely homotopical way by R. Cohen and J. Jones using
a ring spectrum structure on a Thom spectrum of a virtual bundle
over free loop spaces \cite{CJ}. As discovered by S. Voronov
\cite{V}, it comes in fact from a geometric operadic action of the
cacti operad. Very recently J. Klein in \cite{JK} has extended the
homotopy theoretic approach of R. Cohen and J. Jones to Poincar\'e
duality spaces using $A_{\infty}$-ring spectrum technology.
\\
\\
In this paper we adopt a different approach to string topology,
namely we use a geometric version of singular homology introduced
by M. Jakob \cite{MJ1}. And we show how it is possible to define
Gysin morphisms, exterior products and intersection type products
(such as the loop product of M. Chas and D. Sullivan) in the
setting of Hilbert manifolds. Let us point out that three
different types of free loop spaces are used in the mathematical
literature:
\\
- Spaces of continuous loops (\cite{CS} for example),
\\
- Spaces of smooth loops, which are Fr\'echet manifolds but not
Hilbert manifolds (\cite{Br} for some details),
\\
- Spaces of Sobolev class of loops \cite{K} or \cite{GM},
\\
These three spaces are very different from an analytical point of
view, but they are homotopy equivalent. For our purpose, we deal
with Hilbert manifolds in order to have a nice theory of
transversality, the last of maps is the one we use.
\\
\\
In order to perform such intersection theory we recall in section
$2$ what is known about transversality in the context of Hilbert
manifolds. We also describe the manifold structure of free loop
spaces used by W. Klingenberg \cite{K} in order to study closed
geodesics on Riemannian manifolds. The cornerstone of all the
constructions of the next sections will be the "string pull-back",
also used by R. Cohen and J. Jones \cite[diagram 1.1]{CJ}. In
section $2.4$ we extend these techniques $n$-sphere spaces and we
show how to intersect geometrically families of $n$-sphere in $M$.
\\
\\
Section $3$ is devoted to the introduction and main properties
of geometric homology. This theory is based upon bordism classes
of singular manifolds. In this setting families of $n$-sphere in
$M$, which are families parametrized by smooth oriented compact
manifolds, have a clear homological meaning. Of particular interest and crucial importance
for applications to sphere topology is the construction of an explicit Gysin morphism for Hilbert manifolds
in the context of geometric homology (section $3.3$). This construction does not use any
Thom spaces and is based on the construction
of pull-backs for Hilbert manifolds. Such approach seems completely new in this context.
\\
We want to point out that all the constructions performed in this
section work with a generalized homology theory $h_*$ under some
mild assumptions (\cite{MJ1}):
\\
- the associated cohomology theory $h^*$ is multiplicative,
\\
- $h_*$ satisfies the infinite wedge axioms.
\\
\\
In section $4$ the operator $\Delta$, the loop product, the loop
bracket, the intersection morphism and the string bracket are
defined an studied using the techniques introduced in section $2$
and $3$. This section is also concerned with string topology
operations, these operations are parametrized by the topological
space of reduced Sullivan's Chord diagrams $\overline{\mathcal
C\mathcal F^{\mu}_{p,q}(g)}$, which is closely related to the
combinatorics of Riemann surfaces of genus $g$, with $p$-incoming
boundary components and $q$-outgoing. A. Voronov (private
communication) suggested to introduce these spaces of diagrams
because they form a prop and the cacti appear as a sub-operad.
Pushing the work of R. Cohen and V. Godin on the cation of
Sullivan's chord diagrams on free loop spaces further we prove:
\\
\\
{\bf Theorem:} {\it Let $\mathcal LM$ be the free loop space over
a compact $d$ dimensional manifold $M$. For $q>0$ there exist
morphisms:
$$\mu_{n,p,q}(g):H_n(\overline{\mathcal C\mathcal F^{\mu}_{p,q}(g)})\rightarrow
Hom(H_*(\mathcal LM^{\times p}),H_{*+\chi(\Sigma).d+n}(\mathcal
LM^{\times q})).$$ where $\chi(\Sigma)=2-2g-p-q$.}
\\
\\
Moreover as these operations are compatible with the gluing of
reduced Sullivan's diagrams, $H_*(\mathcal LM)$ appear as an
algebra over this Prop. As a corollary one recovers the structure
of a Frobenius algebra on $H_{*+d}(\mathcal LM)$, build in
\cite{CG}, the operator $\Delta$ and the loop product of M. Chas
and D. Sullivan \cite{CS}.
\\
\\
Section $5$ is to devoted to the extension of the results of the
section $4$ to $n$-sphere spaces. In particular, there exists a
commutative and associative product on the homology of these
spaces. The case of $3$-sphere spaces is detailed.
\\
\\
\\
{\bf Acknowledgments:} I would like to thank Andy Baker and Martin
Jakob for their help about infinite dimensional manifolds and
geometric homology. Discussions about string topology with Ralph
Cohen, Yves F{\'e}lix and Jean-Claude Thomas and remarks of Muriel
Livernet were also very useful. I am also thankful to Sasha
Voronov for a careful reading of a preliminary draft of this paper
and for many suggestions. I am very grateful to the algebraic
topology group of Barcelona (CRM, UAB and UB) for organizing a
seminar on this subject, this was my main motivation for writing
this paper. Finally, I warmly thanks the Centre de Recerca de
Matematica for its hospitality.
\section{Infinite dimensional manifolds}

\subsection{Recollections on Hilbert manifolds} This section is
expository, we review the basic facts about Hilbert manifolds, we
refer to \cite{L} (see also \cite{KrM} for a general introduction
to infinite dimensional manifolds).  Moreover all the manifolds we
consider in this paper are Hausdorff and second countable (we need
these conditions in order to consider partitions of unity).

\subsubsection{Differential calculus.} Let $E$ and $F$ be two
topological vector spaces, there is no difficulty to extend the
notion of differentiability of a continuous map between $E$ and
$F$. Hence, let $f:E\rightarrow F$ be a continuous map we say that
$f$ is differentiable at $x\in E$ if for any $v\in E$ the limit:
$$df_x(v)=\lim_{t\rightarrow 0}\frac{f(x+tv)-f(x)}{t}$$
exists. One can define differentials, $C^{\infty}$ morphisms,
diffeomorphisms and so on.

\subsubsection{Hilbert manifolds.} A topological space $X$ is a
manifold modelled on a separable Hilbert space $E$ if there exists
an atlas $\{U_i,\phi_i\}_{i\in I}$ such that:
\\
i) each $U_i$ is an open set of $X$ and $X=\bigcup_{i\in I} U_i$,
\\
ii) $\phi_i:U_i\rightarrow E$ is an homeomorphism,
\\
iii) $\phi_i\phi_j^{-1}$ is a diffeomorphism whenever $U_i\bigcap
U_j$ is not empty.

\subsubsection{Fredholm maps.} A smooth map $f:X\rightarrow
Y$ between two Hilbert manifolds is a Fredholm map if for each
$x\in X$, the linear map
$$df_x:T_x X\longrightarrow T_{f(x)}Y$$ is a
Fredholm operator, that is to say if $ker df_x$ and $coker df_x$
are finite dimensional vector spaces. Recall that the index of a
Fredholm map
$$index:X\rightarrow \mathbb Z$$
$$index (f_x)=dim (ker df_x)-dim (coker df_x)$$
is a continuous map.

\subsubsection{Orientable morphisms.} A smooth map $f:X\rightarrow
Y$ between two Hilbert manifolds is an orientable morphism if it
is a proper map (the pre-image of a compact set is compact) which is also
Fredholm and such that the normal bundle $\nu(f)$ is orientable (for convenience we consider the notion of orientability
with respect to singular homology but we could have chosen to work in a much more general setting).
\\
Let us remark that a closed embedding is an orientable morphism if
and only if $\nu(f)$ is finite dimensional and orientable. A
closed Fredholm map is proper by a result of S. Smale \cite{Sm}.

\subsubsection{Partitions of unity.} A very nice feature of
Hilbert manifolds with respect to Banach manifolds and other type
of infinite dimensional manifolds is the existence of partitions
of unity (\cite[chapter II,3]{L} for a proof). As a consequence
mimicking techniques used in the finite dimensional case, one can
prove that any continuous map
$$f:P\rightarrow X$$
from a finite dimensional manifold $P$ to an Hilbert manifold $X$
is homotopic to a $C^{\infty}$ one. And we can also smoothen
homotopies.

\subsection{Transversality} We follow the techniques developed by A. Baker and C. {\"O}zel in
\cite{BO} in order to deal with transversality in an infinite
dimensional context.

\subsubsection{Transversal maps} Let $f:X\rightarrow
Y$ and $g:Z\rightarrow Y$ be smooth maps between two Hilbert
manifolds. Then they are transverse at $y\in Y$ if

$$df(T_x X)+dg(T_z Z)=T_y Y$$
with $f(x)=g(z)=y$. The maps are transverse if they are tranverse
at any point $y\in Imf\cap Img$.

\subsubsection{Pull-backs} Let us recall the main results about
pull-backs of Hilbert manifolds. We consider the following
diagram:

\begin{diagram}
Z &\lTo^{g^*f}&  Z\cap_{Y} X\\
\dTo^{g}&&\dTo_{\phi}\\
Y&\lTo_{f}& X\\
\end{diagram}

\noindent where $Z$ is a finite dimensional manifold and
$f:X\rightarrow Y$ is an orientable map.
\\
\\
Using an infinite dimensional version of the implicit function
theorem \cite[Chapter I,5]{L}, one can prove the following result:

\subsubsection{\bf Proposition}{\it \cite[prop. 1.17]{BO} If the map
$$f:X\rightarrow Y$$ is an
orientable morphism and
$$g:Z\rightarrow Y$$ is a smooth map
transverse to $f$, then the pull-back map:
$$g^*f:Z\cap_{Y} X\longrightarrow Z$$
is an orientable morphism.}
\\
\\
Moreover the finite dimensional hypothesis on $Z$ enables to
prove:

\subsubsection{\bf Theorem} {\it \cite[Th. 2.1, 2.4]{BO} Let
$$f:X\rightarrow Y$$ be an orientable
morphism and let $$g:Z\rightarrow Y$$ be a smooth map from a
finite dimensional manifold $Z$. Then $g$ can be deformed by a
smooth homotopy until it is transverse to $f$.}

\subsection{Free loop spaces}
If we want to do intersection theory with spaces of closed curves,
we need to consider them as smooth manifolds. Following
\cite[Chapter 3]{Br}, one can consider the space
$C^{\infty}(S^1,M)$ of all smooth curves as an Inverse Limit
Hilbert manifold. But we prefer to enlarge this space and to
consider ${\bf H}^1(S^1,M)$ the space of ${\bf H}^1$ curves. This space has the advantage
to be an Hilbert manifold. With this choice we can apply all the techniques
described in the sections $2.1$ and $2.2$.
\\
In fact, we could also consider ${\bf H}^n(S^n,M)$ the space of ${\bf H}^n$-curves, all these spaces are
Hilbert manifolds, and they are also all homotopy equivalent to the ILH manifold $C^{\infty}(S^1,M)$.
And all these manifolds are homotopy equivalent to the space of continuous maps $C^0(S^1,M)$ equipped with the
compact open topology.

\subsubsection{Manifold structure} In order to define an Hilbert manifold structure
on free loop spaces we follow
W. Klingenberg's approach \cite{K}.
\\
Let $M$ be a simply-connected Riemannian manifold of dimension
$d$. We set
$$\mathcal LM={\bf H}^1(S^1,M).$$
The manifold $\mathcal LM$ is formed by the continuous curves
$\gamma:S^1\rightarrow M$ of class ${\bf H}^1$, it is modelled on the
Hilbert space $\mathcal L\mathbb R^d={\bf H}^1(S^1,\mathbb R^d)$. The
space $\mathcal L\mathbb R^d$ can be viewed as the completion of
the space $C^{\infty}_{p}(S^1,\mathbb R^d)$ of piecewise
differentiable curves with respect to the norm $\|-\|_1$. This norm is defined
via the scalar product:
$$<\gamma, \gamma'>_1=\int\gamma(t)\diamond\gamma'(t)dt+\int \delta\gamma(t)\diamond \delta\gamma'(t) dt,$$
where $\diamond$ is the canonical scalar product of $\mathbb R^d$. As $S^1$
is $1$-dimensional, we notice that by Sobolev's embedding theorem
elements of $\mathcal L\mathbb R^d$ can be represented by continuous curves.
\\
Let us describe an atlas $\{U_{\gamma},\phi_{\gamma}\}$ of $\mathcal LM$. We
take $\gamma \in C^{\infty}_{p}(S^1,M)$ a piecewise differentiable
curve in $M$ (notice that $C^{\infty}_{p}(S^1,M)\subset
{\bf H}^1(S^1,M)$) and consider the pullback:
\begin{diagram}
\gamma^*TM&\rTo&TM\\
\dTo&&\dTo\\
S^1&\rTo_{\gamma}&M.\\
\end{diagram}
Now let $T^{\epsilon}_{\gamma}\subset \gamma^*TM$ be the set of vectors of
norm less than $\epsilon$. The exponential map
$$exp: T^{\epsilon}_{\gamma}\longrightarrow M$$
identifies $T^{\epsilon}_{\gamma}$ with an open set of $M$ and induces a
map:
$$\mathcal LT^{\epsilon}_{\gamma}\longrightarrow\mathcal LM.$$
moreover as $\gamma^*TM$ is a trivial vector bundle (because M is
$1$-connected), we fix a trivialization:
$$\varphi:\gamma^*TM\longrightarrow S^1\times \mathbb{R}^d$$
this gives a chart:
$$\phi_{\gamma}:\mathcal LT^{\epsilon}_{\gamma}\longrightarrow\mathcal L\mathbb{R}^d.$$
\subsubsection{\bf Remark} In fact, the manifold structure on $\mathcal
LM$ does not depend on a choice of a particular Riemannian metric
on $M$.
\subsubsection{The tangent bundle.} Let $TM\rightarrow M$ be the
tangent bundle of $M$. The tangent bundle of $\mathcal LM$ denoted
by $T\mathcal LM$ can be identified with $\mathcal LTM$, this is
an infinite dimensional vector bundle where each fiber is
isomorphic to the Hilbert space $\mathcal L\mathbb R^d$. Let
$\gamma\in \mathcal LM$ we have:
$$T\mathcal LM_{\gamma}=\Gamma(\gamma^*TM),$$
where $\Gamma(\gamma^*TM)$ is the space of sections of the
pullback of the tangent bundle of $M$ along $\gamma$ (this is the
space of ${\bf H}^1$ vector fields along the curve $\gamma$). A
trivialization $\varphi$ of $\gamma^*TM$ induces an isomorphism:
$$T\mathcal LM_{\gamma}\cong \mathcal L\mathbb R^d.$$
The tangent bundle of $\mathcal LM$ has been studied in \cite{CSt}
and \cite{Mo}.

\subsubsection{Riemannian structure.} The manifold $\mathcal LM$
has a natural Riemannian metric, the scalar product on $T\mathcal
LM_{\gamma}\cong \mathcal L\mathbb R^d$ comes from $<-,->_1$.

\subsubsection{The $S^1$-action.} The circle acts on $\mathcal
LM$:
$$\Theta:S^1\times \mathcal LM\longrightarrow \mathcal LM$$
by reparametrization:
$$\Theta(\theta,\gamma):t\mapsto \gamma(t+\theta).$$
Of course this action is not free, it is continuous but not
differentiable.
\\
Let $\gamma\in \mathcal LM$, then the isotropy subgroup
$Is(\gamma)$ of $\gamma$ is $S^1$ if and only if $\gamma$ is a
constant map, in this case we say that $\gamma$ is of multiplicity
$0$. Otherwise it is isomorphic to a finite cyclic group and the
multiplicity of the curve is equal to the order of the isotropy
subgroup. Let $\mathcal LM^{(m)}$ be the space of curves of
multiplicity equal to $m$. This gives an $S^1$-equivariant
partition of $\mathcal LM$:
$$\mathcal LM=\bigcup_{m} \mathcal LM^{(m)}.$$
The space $\mathcal LM^{(0)}$ can be identified with $M$, and
$\mathcal LM^{(1)}$ is called the space of prime curves.

\subsubsection{The string pullback.} Let us consider the evaluation
map
$$ev_0:\mathcal LM\rightarrow M$$
$$\gamma\mapsto\gamma(0),$$
this is a submersion of Hilbert manifolds (this follows
immediately from the definition of $T\mathcal LM$). As the map
$ev_0\times ev_0$ is transverse to the diagonal map $\Delta$
(because $ev_0\times ev_0$ is a submersion), we can form the {\it
string pull-back \cite[(1.1)]{CJ}}:

\begin{diagram}
\mathcal LM\times \mathcal LM &\lTo^{\widetilde{\Delta}}& \mathcal
LM\cap_{M}\mathcal LM\\
\dTo^{ev_0\times ev_0}&&\dTo_{ev}\\
M\times M&\lTo^{\Delta}&M,\\
\end{diagram}

\noindent by transversality this is a diagram of Hilbert
manifolds. We have:
$$\mathcal
LM\cap_{M}\mathcal LM=\{(\alpha,\beta)\in  \mathcal LM\times
\mathcal LM / \alpha(0)=\beta(0) \}.$$ The map
$$\widetilde{\Delta}:\mathcal
LM\cap_{M}\mathcal LM\rightarrow \mathcal LM\times \mathcal LM$$
is a closed embedding of codimension $d$.
\\
As the normal bundle $\nu_{\tilde \Delta}$ is the pull-back of $\nu_{\Delta}$ and as this last one is
isomorphic to
$TM$ , we
deduce that $\tilde \Delta$ is an orientable morphism.

\subsubsection{Families of closed strings} A family of closed
strings in $M$ is a smooth map
$$f:P\rightarrow \mathcal LM$$
from a compact orientable manifold $P$.
\\
\\
The proposition below gives a necessary but not sufficient
condition in order to do intersection of families of closed
strings.

\subsubsection{\bf Proposition} {\it If $P\times Q \stackrel {f\times g}{\rightarrow} \mathcal LM\times \mathcal LM$
is transverse to $\widetilde{\Delta}$ then $ev_0f$ and $ev_0g$ are
transverse in $M$.}
\\
\\
Now we suppose that $(P,f)$ and $(Q,g)$ are two orientable compact
manifolds of dimensions $p$ and $q$ respectively. Moreover we
suppose that they are such that $f\times g$ is transverse to
$\widetilde{\Delta}$. We denote by $P\ast Q$ the pullback:
\begin{diagram}
P\times Q &\lTo& P\ast Q\\
\dTo^{f\times g}&& \dTo_{\psi}\\
 \mathcal LM\times \mathcal LM
&\lTo^{\widetilde{\Delta}}& \mathcal
LM\cap_{M}\mathcal LM.\\
\end{diagram}
Then $P\ast Q$ is a compact orientable submanifold of $P\times Q$
of dimension $p+q-d$.

\subsubsection{Composing loops} Let us define the map:
$$\Upsilon:\mathcal LM\cap_{M}\mathcal LM\longrightarrow \mathcal LM.$$
Let $(\alpha,\beta)$ be an element of $\mathcal LM\cap_{M}\mathcal
LM$ then $\Upsilon(\alpha,\beta)$ is the curve defined by:
\\
$\Upsilon(\alpha,\beta)(t)=\alpha(2t)$ if $t\in [0,1/2]$
\\
$\Upsilon(\alpha,\beta)(t)=\beta(2t-1)$ if $t\in [1/2,1].$
\\
We notice that this map is well defined because we compose
piecewise differential curves, hence no "dampening" constructions
are needed as in \cite[remark about construction (1.2)]{CJ}.
\\
The construction of $\Upsilon$ comes from the co-H-space structure
of $S^1$ i.e. the pinching map:
$$S^1\longrightarrow S^1\vee S^1.$$

\subsubsection{Intersection of families of closed strings} Now
consider two families of closed strings $(P,f)$ and $(Q,g)$, by
deforming $f\times g$ one can produce a new family of closed
strings $(P\ast Q,\Upsilon\psi)$ in $M$. We also notice that the
image of $\Upsilon\psi$ lies in $\mathcal LM^{(0)}\bigcup \mathcal
LM^{(1)}$.

\subsubsection{\bf Remark} All we have done with free loop spaces can
be performed for manifolds of maps from a space which is a co-H-space and
a compact orientable manifold to a compact Riemannian manifold.

\subsection{$n$-sphere spaces.} Let $M$ be a $n$-connected $d$-dimensional
compact oriented smooth manifold.

\subsubsection{\bf Definition} {\it We call the $n$-sphere space
of $M$ and we denote it by $\mathcal S_{n}M$ the space of
${\bf H}^n$-maps from $S^n$ to $M$.}

\subsubsection{\bf Remark} By Sobolev's embedding theorem we know
that $${\bf H}^n(S^n,M)\subset C^0(S^n,M).$$

\subsubsection{\bf Theorem} {\it $\mathcal S_{n}M$ is an Hilbert manifold.}
\\
\\
{\bf Proof} As for free loop spaces an atlas of $\mathcal S_{n}M$
can be given by considering $\bf H^n$ vector fields along all maps
$$b:S^n\rightarrow M.$$
Then using a trivialization of $b^*TM$ we deduce that $\mathcal
S_{n}M$ is modelled on the separable Hilbert space $\mathcal
S_{n}\mathbb R^d$. $\square$
\\
\\
The tangent bundle of $\mathcal S_{n}M$ has the same description
as $T\mathcal LM$. It can be identified with $\mathcal S_{n}TM$
and we have: $$T\mathcal S_{n}M_{b}=\Gamma(b^*TM).$$ Moreover,
$\mathcal S_{n}M$ is a Riemannian manifold.

\subsubsection{The $n$-sphere pull-back.} Let fix a base point $0$
in $S^n$, the evaluation map:
$$ev_0:\mathcal S_{n}M\rightarrow M$$
is clearly a submersion of Hilbert manifolds, then we can form the
pull-back of Hilbert manifolds:
\begin{diagram}
\mathcal S_{n}M\times \mathcal S_{n}M &\lTo^{\widetilde{\Delta}}&
\mathcal
S_{n}M\cap_{M}\mathcal S_{n}M\\
\dTo^{ev_0\times ev_0}&&\dTo_{ev}\\
M\times M&\lTo^{\Delta}&M,\\
\end{diagram}
and the map $\tilde \Delta$ is an orientable morphism.

\subsubsection{Composing $n$-sphere} Thanks to the pinching map:
$$S^n\longrightarrow S^n\vee S^n$$
one can define:
$$\Upsilon:S_{n}M\cap_{M}\mathcal S_{n}M\longrightarrow\mathcal S_{n}M.$$

\subsubsection{Intersection of families of $n$-sphere} As for families of closed strings,
we consider two families of $n$-sphere in $M$ denoted by $(P,f)$
and $(Q,g)$, by deforming $f\times g$ and taking the pullback
$P\ast Q$ along $\tilde \Delta$ one can produce a new family of
$n$-sphere $(P\ast Q,\Upsilon\psi)$ in $M$.

\section{Geometric homology theories}

As R. Thom proved it is not possible in general to represent
singular homology classes of a topological space $X$ by singular
maps i.e continuous maps:
$$f:P\longrightarrow X$$
from an oriented manifold to $X$. But, M. Jakob in \cite{MJ1},
\cite{MJ2} proves that if we add a cohomological information to
the map $f$ (a singular cohomological class of $P$), then
Steenrod's realizability problem with this additional cohomological data has an affirmative answer. In
these two papers he develops a geometric version of homology. This
geometric version seems to be very nice to deal with Gysin
morphisms, intersection products and so on.
\\
All the constructions we give below and also their applications to string topology work out for more general
homology theories: bordism, topological $K$-theory for example. We refer the reader to
\cite{MJ1}, \cite{MJ2} and \cite{MJ3} for the definitions of these geometric theories.

\subsection{An alternative description of singular homology}

\subsubsection{Geometric cycles} Let $X$ be a topological space, a
geometric cycle is a triple $(P,a,f)$ where:
$$f:P\longrightarrow X$$
is a continuous map from a compact connected orientable manifold
$P$ to $X$ (i.e a singular manifold over $X$), and $a\in
H^*(P,\mathbb Z)$. If $P$ is of dimension $p$ and $a\in
H^m(P,\mathbb Z)$ then $(P,a,f)$ is a geometric cycle of degree
$p-m$. Take the free abelian group generated by all the geometric
cycles and impose the following relation:
$$(P,\lambda.a+\mu.b,f)=\lambda.(P,a,f)+\mu.(P,b,f).$$
Thus we get a graded abelian group.

\subsubsection{Relations} In order to recover singular homology we
must impose the two following relations on geometric cycles:
\\
\\
i) ({\bf Bordism relation}) If we have a map $h:W\rightarrow X$
where $W$ is an orientable bordism between $(P,f)$ and $(Q,g)$
i.e. $$\partial W=P\bigcup Q^{-}.$$ Let $i_1:P\hookrightarrow W$
and $i_2:Q\hookrightarrow W$ be the canonical inclusions, then for
any $c\in H^*(W,\mathbb Z)$ we impose:
$$(P,i_1^*(c),f)=(Q,i_2^*(c),g).$$
\\
ii) ({\bf Vector bundle modification}) Let $(P,a,f)$ be a
geometric cycle and consider a smooth orientable vector bundle
$E\stackrel{\pi}{\rightarrow}P$, take the unit sphere bundle
$S(E\oplus 1)$ of the Whitney sum of $E$ with a copy of the
trivial bundle over $M$. The bundle $S(E\oplus 1)$ admits a
section $\sigma$, by $\sigma_{!}$ we denote the Gysin morphism in
cohomology associated to this section. Then we impose:
$$(P,a,f)=(S(E\oplus 1),\sigma_{!}(a),f\pi).$$
\\
An equivalence class of geometric cycle is denoted by $[P,a,f]$,
let call it a geometric class. And $H'_{q}(X)$ is the abelian
group of geometric classes of degree $q$.

\subsubsection{\bf Theorem}{\it \cite[Cor. 2.36]{MJ1} The morphism:
$$H'_{q}(X)\longrightarrow H_q(X,\mathbb Z)$$
$$[P,a,f]\mapsto f_*(a\cap [P])$$
where $[P]$ is the fundamental class of $P$ is an isomorphism of
abelian groups.}

\subsection{Cap product and Poincar\'e duality \cite[3.2]{MJ2}} The
cap product between $H^*(X,\mathbb Z)$ and $H'_*(X)$ is given by
the following formula:
$$\cap: H^p(X,\mathbb Z)\otimes H'_q(X)\longrightarrow H'_{q-p}(X)$$
$$u\cap [P,a,f]=[P,f^*(u)\cup a,f].$$
Let $M$ be a $d$-dimensional smooth compact orientable manifold
without boundary then the morphism:
$$H^p(M,\mathbb Z)\longrightarrow H'_{d-p}(M)$$
$$x\mapsto [M,x,Id_M]$$
is an isomorphism.

\subsection{Gysin morphisms}(\cite{MJ3} for a finite dimensional version) We want to consider Gysin
morphisms in the context of infinite dimensional manifolds.
\\
Let us recall two possible definitions for Gysin morphisms in the finite dimensional context. The following one is
only relevant to the final dimensional case. Let us take a morphism:
$$f:M^m\longrightarrow N^n$$
of Poincar\'e duality spaces. Then we define:
$$f_{!}:H_*(N^n)\stackrel{D}{\rightarrow}H^{n-*}(N^n)\stackrel{f}{\rightarrow}H^{n-*}(M^n)
\stackrel{D^{-1}}{\rightarrow}H_{*+m-n}(M^m),$$ where $D$ is the
Poncar\'e duality isomorphism.
\\
For the second construction, if $f$ is an embedding of smooth oriented manifold
then one can apply the Pontryagin-Thom collapse $c$ to the Thom space of the normal bundle of $f$ and then apply
the Thom isomorphism $th$:
$$f_{!}:H_*(N^n)\stackrel{c}{\rightarrow}H_{*}(Th(\nu(f)))
\stackrel{th}{\rightarrow}H_{*+m-n}(M^m)$$ In the infinite
dimensional context, we can not use Poincar\'e duality, and to be
as explicit as possible we do not want to use the Pontryagin-Thom
collapse and the Thom isomorphism. We prefer to use a very
geometrical interpretation of the Gysin morphism which is to take
pull backs of cycles along the map $f$.
\\
\\
So, we
take $i:X\rightarrow Y$ an orientable morphism of Hilbert
manifolds and we suppose that $\nu(i)$ is $d$-dimensional. Let us
define:
$$i^{!}:H'_p(Y)\longrightarrow H'_{p-d}(X).$$
Let $[P,a,f]$ be a geometric class in $H'_{p}(Y)$, we can choose a
representing cycle $(P,a,f)$. If $f$ is not smooth, we know that
it is homotopic to a smooth map by the existence of partitions of
unity on $Y$, moreover we can choose it transverse to $i$, by the
bordism relation all these cycles represent the same class. Now we
can form the pull-back:
\begin{diagram}
P &\lTo^{f^*i}&  P\cap_{Y} X\\
\dTo^{f}&&\dTo_{\phi}\\
Y&\lTo_{i}& X\\
\end{diagram}
we set:
$$i^{!}([P,a,f])=(-1)^{d.|a|}[P\cap_{Y} X,(f^*i)^*(a),\phi].$$
The sign is taken from \cite[3.2c)]{MJ3}, the Gysin morphism can be viewed as a product  for
bivariant theories \cite{FM}.

\subsection{The cross product \cite[3.1]{MJ2}} The cross product is
given by the pairing:
$$\times:H'_q(X)\otimes H'_p(Y)\longrightarrow H'_{p+q}(X\times Y)$$
$$[P,a,f]\times [Q,b,g]=(-1)^{dim(P).|b|}[P\times Q,a\times b, f\times g].$$
The sign makes the cross product commutative. let
$$\tau:X\times Y\rightarrow Y\times X$$
be the interchanging morphism then:
$$\tau_*(\alpha\times \beta)=(-1)^{|\alpha||\beta|}\beta\times\alpha.$$

\subsection{The intersection product(\cite[sect.3]{MJ3})} Let us return to the finite dimensional case and consider $M$
an orientable compact $d$-dimensional manifold. Like for Gysin
morphisms in order to be very explicit we avoid the classical
constructions of the intersection product that use either
Poincar\'e duality or the Thom isomorphism.
\\
Let $[P,x,f]\in
H'_{n_1}(M)$ and $[Q,y,g]\in H'_{n_2}(M)$, we suppose that $f$ and
$g$ are transverse in $M$, then we form the pull back:
\begin{diagram}
P\times Q &\lTo^{j}&  P\cap_{M} Q\\
\dTo^{f\times g}&&\dTo_{\phi}\\
M\times M&\lTo_{\Delta}& M\\
\end{diagram}
and define the pairing:
$$-\bullet-:H'_{n_1}(M)\otimes H'_{n_2}(M)\stackrel{\times}{\rightarrow}H'_{n_1+n_2}(M\times M)
\stackrel{\Delta^{!}}{\rightarrow} H'_{n_1+n_2-d}(M).$$ Hence, we
set:
$$[P,a,f]\bullet[Q,b,g]=(-1)^{d.(|a|+|b|)+ dim(P).|b|}[P\cap_{M} Q,j^*(a\times b),\phi].$$
Let $l:P\cap_{M} Q\rightarrow P$ and $r:P\cap_{M} Q\rightarrow Q$
be the canonical maps, then we also have:
$$[P,a,f]\bullet[Q,b,g]=(-1)^{d.(|a|+|b|)+dim(P).|b|}[P\cap_{M} Q,l^*(a)\cup r^*(b),\phi].$$
With these signs conventions the intersection product $\bullet$
makes $H'_{*+d}(M)$ into a graded commutative algebra:
$$[P\cap_{M} Q,l^*(a)\cup r^*(b),\phi]=(-1)^{(d-dim(P)-|a|)(d-dim(Q)-|b|)}[Q\cap_{M} P,l^*(b)\cup r^*(a),\phi].$$

\section{String topology}
In this section, using the theory of geometric cycles we show how
to recover the $BV$-structure on
$$\mathbb H_*(\mathcal LM):= H'_{*+d}(\mathcal LM,\mathbb Z)$$
introduced in \cite{CS} and studied from a homotopical point of
view in \cite{CJ}.
\\
We also define the intersection morphism, the string bracket of \cite{CS} and
string topology operations (we extend the Frobenius structure
given in \cite{CG} to a homological action of the space of Sullivan's chord diagrams).
\\
\\
{\bf Remark:} In this section we use the language of operads and
algebras over an operad (in order to state some results in a nice
and appropriate framework). For definitions and examples of
operads and algebras over an operad we refer to \cite{GJ},
\cite{GK}, \cite{KM}, \cite{Lo}, \cite{MSS} and \cite{V}.

\subsection{The operator $\Delta$} First we define the
$\Delta$-operator on $\mathbb H_*(\mathcal LM)$. Let us consider a
geometric cycle $[P,a,f]\in H'_{n+d}(\mathcal LM)$, we have a map:
$$ \Theta_f:S^1\times P\stackrel{Id\times f}{\rightarrow}S^1\times \mathcal LM
\stackrel{\Theta}{\rightarrow} \mathcal LM.$$

\subsubsection{\bf Definition} {\it There is an operator
$$\Delta:H'_{n+d}(\mathcal LM)\rightarrow H'_{n+d+1}(\mathcal LM)$$
given by the following formula:
$$\Delta([P,a,f])=(-1)^{|a|}[S^1\times P,1\times a,\Theta_f].$$}

\subsubsection{\bf Proposition \cite[prop. 5.1]{CS}} {\it The operator verifies:
$\Delta^2=0$.}
\\
\\
{\bf Proof} This follows from the associativity of the cross
product and the nullity of $[S^1\times S^1,1\times 1,\mu]\in
H_2'(S^1)$ where $\mu$ is the product on $S^1$. $\square$

\subsection{Loop product} Let us take $[P,a,f]\in
H'_{n_1+d}(\mathcal LM)$ and $[Q,b,g]\in H'_{n_2+d}(\mathcal LM)$.
We can smooth $f$ and $g$ and make them transverse to
$\widetilde{\Delta}$, then we form the pull-back $P\ast Q$.

\subsubsection{\bf Definition}{\it
Let $l:P\ast Q\rightarrow P$ and $r:P\ast Q\rightarrow Q$ be the
canonical maps, then we have the pairing:
$$-\bullet-:H'_{n_1+d}(\mathcal LM)\otimes H'_{n_2+d}(\mathcal LM)\longrightarrow
H'_{n_1+n_2+d}(\mathcal LM)$$
$$[P,a,f]\bullet[Q,b,g]=(-1)^{d.(|a|+|b|)+dim(P).|b|}[P\ast Q,l^*(a)\cup r^*(b),\Upsilon
\psi],$$ let call it the loop product.}

\subsubsection{\bf Proposition \cite[Thm. 3.3]{CS}}{\it The loop product is associative
and commutative.}
\\
\\
{\bf Proof} The associativity of the loop product follows from the
associativity of the intersection product, the cup product and the
fact that $\Upsilon$ is also associative up to homotopy.
\\
In order to prove the commutativity of $\bullet$ we follow the
strategy of \cite[Lemma 3.2]{CS}.
\\
There is a smooth interchanging map:
$$\tau:\mathcal LM\cap_{M}\mathcal LM \rightarrow \mathcal LM\cap_{M}\mathcal LM.$$
Let $[P,a,f]$ and $[Q,b,g]$ be two geometric classes the formula
of \cite[lemma 3.2]{CS} gives an homotopy $H$ (so this is also a
bordism) between:
$$P\ast Q\stackrel{\psi}{\rightarrow}\mathcal LM\cap_{M}\mathcal LM\stackrel{\Upsilon}{\rightarrow}\mathcal LM$$
and
$$P\ast Q\stackrel{\psi}{\rightarrow}\mathcal LM\cap_{M}\mathcal LM\stackrel{\Upsilon\tau}{\rightarrow}\mathcal LM.$$

This bordism identifies $[P\ast Q,l^*(a)\cup r^*(b),\Upsilon
\psi]$ and $$[P\ast Q,\tau^*(l^*(a)\cup r^*(b)),\Upsilon\tau
\psi]$$ which is equal to:
$$(-1)^{(dim(P)-d-a)(dim(P)-d-b)}[Q\ast P,l^*(b)\cup
r^*(a),\Upsilon\psi].$$
 $\square$

\subsection{Loop bracket} Let $[P,a,f]$ and $[Q,b,g]$ be two geometric
classes, in the preceding section we have defined a bordism
between
$$P\ast Q\stackrel{\psi}{\rightarrow}\mathcal LM\cap_{M}\mathcal LM\stackrel{\Upsilon}{\rightarrow}\mathcal LM$$
and
$$P\ast Q\stackrel{\psi}{\rightarrow}\mathcal LM\cap_{M}\mathcal LM\stackrel{\Upsilon\tau}{\rightarrow}\mathcal LM.$$
Using the same homotopy one can define another bordism between
$$P\ast Q\stackrel{\psi}{\rightarrow}\mathcal LM\cap_{M}\mathcal LM\stackrel{\Upsilon\tau}{\rightarrow}\mathcal LM$$
and
$$P\ast Q\stackrel{\psi}{\rightarrow}\mathcal LM\cap_{M}\mathcal LM\stackrel{\Upsilon\tau^2}{\rightarrow}\mathcal LM.$$
Composing these two bordisms one obtains a geometric class:
$$(-1)^{|a|+|b|}[S^1\times P\ast Q,1\times l^*(a)\cup r^*(b)\times 1, \tilde H].$$

\subsubsection{\bf Definition}{\it The loop bracket is the pairing:
$$\{-,-\}:H'_{n_1+d}(\mathcal LM)\otimes H'_{n_2+d}(\mathcal LM)\longrightarrow
H'_{n_1+n_2+d+1}(\mathcal LM)$$
$$\{[P,a,f],[Q,b,g]\}=(-1)^{(d+1).(|a|+|b|)+dim(P).|b|}[S^1\times P\ast Q,1\times l^*(a)\cup r^*(b),\tilde H].$$ }
\\
\\
There is another way to define the bracket by setting \cite[Cor.
5.3]{CS}:
$$\{\alpha,\beta\}=(-1)^{|\alpha|}\Delta(\alpha\bullet\beta)-(-1)^{|\alpha|}\Delta(\alpha)\bullet\beta-\alpha\bullet
\Delta(\beta).$$
Together with this bracket, $(\mathbb H_*(\mathcal LM),\bullet,\{-,-\})$ is a Gerstenhaber algebra.

 \subsubsection{\bf Theorem}{\cite[Thm. 4.7]{CS}}{\it The triple $(\mathbb H_*(\mathcal LM),\bullet,\{-,-\})$
is a Gerstenhaber algebra:
\\
i) $(\mathbb H_*(\mathcal LM),\bullet)$ is a graded associative and commutative algebra.
\\
ii) The loop bracket $\{-,-\}$ is a Lie bracket of degree $+1$:
$$\{\alpha,\beta\}=(-1)^{(|\alpha|+1)(|\beta|+1)}\{\beta,\alpha\},$$
$$\{\alpha,\{\beta,\gamma\}\}=\{\{\alpha,\beta\}.\gamma\}+(-1)^{(|\alpha|+1)(|\beta|+1)}
\{\beta,\{\alpha,\gamma\}\},$$
\\
iii) $\{\alpha,\beta\bullet\gamma\}=\{\alpha,\beta\}\bullet \gamma+(-1)^{|\beta|(|\alpha|+1)}
\beta\bullet\{\alpha,\gamma\}.$}

\subsubsection{\bf Remark.} Let us recall that there are two important examples of Gerstenhaber algebras:
\\
- The first one is the
Hochschild cohomology of a differential graded associative algebra $A$:
$$HH^*(A,A),$$
this goes back to M. Gerstenhaber \cite{Ge}.
\\
- The second example is the singular homology of a double loop space:
$$H_*(\Omega^2X),$$
this is due to F. Cohen \cite{CLM}.
\\
Both examples are related by the Deligne's conjecture proved in
many different ways by C. Berger and B. Fresse \cite{BF}, M.
Kontsevich and Y. Soibelman \cite{Ko2}, J. McClure and J. Smith
\cite{MS}, D. Tamarkin \cite{T} and S. Voronov \cite{Vd} (see also
M. Kontsevich \cite{Ko1}). This conjecture states that there is a
natural action of an operad $C_2$ quasi-isomorphic to the chain
operad of little $2$-discs on the Hochschild cochain complex of an
associative algebra.
\\
Hochschild homology enters the theory by the following results of R. Cohen and J.D.S. Jones \cite[Thm. 13]{CJ}:
\\
if $C^*(M)$ denotes the singular cochains of a manifold $M$, then there is an isomorphism of associative algebras:
$$HH^*(C^*(M),C^*(M))\cong \mathbb H_{*}(\mathcal LM).$$

\subsection{The $BV$-structure} In \cite{CJ} and \cite{CS} it is proved that $\mathbb
H_{*}(\mathcal LM)$ is a $BV$-algebra (we refer to \cite{G} for
$BV$-structures).

\subsubsection{\bf Theorem \cite[Th. 5.4]{CS}}{\it The loop product $\bullet$ and the operator
$\Delta$ makes $\mathbb H_{*}(\mathcal LM)$ into a Batalin
Vilkovisky algebra, we have the following relations:
\\
i) $(\mathbb H_{*}(M^{S^1}),\bullet),$ is a graded commutative
associative algebra.
\\
ii) $\Delta^2=0$
\\
iii)
$(-1)^{|\alpha|}\Delta(\alpha\bullet\beta)-(-1)^{|\alpha|}\Delta(\alpha)\bullet\beta-\alpha\bullet
\Delta(\beta)$ is a derivation of each variable.}
\\
\\
The proof of the theorem in the context of geometric homology is
given by building explicit bordisms between geometric cycles. All
these bordisms are described in \cite{CS}.

\subsubsection{\bf Remark.}E. Getzler introduced $BV$-algebras in the context of $2$-dimensional
topological field theories \cite{G}.
And he proved that $H_*(\Omega^2M)$ is a $BV$-algebra if $M$ has a $S^1$ action.
Other examples are provided by the de Rham cohomology
of manifolds with $S^1$-action.
\\
The $BV$-structure on $\mathbb H_*(\mathcal LM)$ comes in fact
form a geometric action of the $cacti$ operad \cite{CJ}, \cite{V}
(normalized cacti with spines in the terminology of R. Kaufmann
\cite{K}). Roughly speaking an element of $cacti(n)$ is a
tree-like configuration of $n$-marked circles in the plane. The
cacti operad is homotopy equivalent to the little framed discs
operad \cite{V}. And we know since the work of E. Getzler that the
homology of the little framed discs operad gives the $BV$ operad
\cite{G}.
\\
Let us explain this geometric action. First let us define the
space $\mathcal L^{cacti(n)}M$ (denoted by $L_kM$ in \cite{CJ})
as:
$$\mathcal L^{cacti(n)}M=\{(c,f):c\in cacti(n), f:c\rightarrow M\}$$
we take the Gysin morphism along a map:
$$cacti(n)\times \mathcal LM^{\times n}\longleftarrow \mathcal L^{cacti(n)}M$$
to any element $c\in cacti(n)$ one can associate a map:
$$S^1\rightarrow c$$
then we get:
$$\mathcal L^{cacti(n)}M\longrightarrow \mathcal LM.$$
For $n=1,2$ we know from R. Kaufmann's description of $cacti$
\cite{RK} that $cacti(n)$ is a smooth manifold. In that case all
the maps defined above are maps of Hilbert manifolds and they give
also a very nice description of the action of $H'_*(cacti)$ on
$H'_{*+d}(\mathcal LM)$.
\\
So, it is certainly worth building a smooth structure on $cacti$
or on an operad homotopy equivalent that acts in the same way.
This would give a more conceptual proof of the preceding theorem.

\subsection{Constant strings} We have a canonical embedding:
$$c:M\hookrightarrow \mathcal LM$$
$c$ induces a map:
$$c_*:H'_{n+d}(M)\rightarrow H'_{n+d}(\mathcal LM).$$
The morphism $c_*$ is clearly a morphism of commutative algebras.

\subsection{Intersection morphism} Let recall that the map
$$ev_0:\mathcal LM\longrightarrow M$$
is a submersion (in fact this is a smooth fiber bundle of Hilbert
manifolds). Hence if we choose a base point $m\in M$ the fiber of
$ev_0$ in $m$ is the Hilbert manifold $\Omega M$ of based loops in
$M$. Consider the morphism:
$$i:\Omega M\hookrightarrow \mathcal LM$$
from the based loops in $M$ to the free loops in $M$,
this is an orientable morphism of codimension $d$.
\\
Let us describe the intersection morphism:
$$I=i^{!}:\mathbb
H_*(\mathcal LM)\rightarrow H_*(\Omega M).$$ Let $[P,a,f]\in
H'_{n+d}(\mathcal LM)$ be a geometric class, one can define
$I([P,a,f])$ in two ways:
\\
\\
i) using the Gysin morphism : $I([P,a,f])=(-1)^{d.|a|}[P\cap_{\mathcal
LM}\Omega M,(f^*i)^*(a),\phi]$.
\\
\\
A better way is certainly to notice that this is the same as doing
the loop product with $[c_m,1,c]$ where $c_m$ is a point and
$c:c_m\rightarrow \mathcal LM$ is the constant loop space at the
point $m$, then we have:
\\
\\
ii) $I([P,a,f])=(-1)^{d.|a|}[P\ast c_m,l^*(a),\psi]$.
\\
\\
We remark that $P\ast c_m$ is either empty (depending on the
dimension of $P$, for example when $dim P<d$) or equal to $m$. And
if $|a|>0$ we also have $I([P,a,f])=0$.

\subsubsection{\bf Proposition \cite[Prop 3.4]{CS}}{\it The intersection morphism $I$ is a morphism of
associative algebras.}
\\
\\
{\bf Proof.} The algebra structure on $H'_*(\Omega M)$ comes from
the Pontryagin product which is the restriction of $\Upsilon$ to
$\Omega M\times \Omega M$, we have the following diagram:
\begin{diagram}
\Omega M\times\Omega M &\rTo^{\Upsilon_{\Omega M\times\Omega M}}& \Omega M\\
\dTo^{i\times i}&&\dTo^{i}\\
\mathcal LM\cap_{M}\mathcal LM&\rTo^{\Upsilon}&\mathcal LM.\\
\end{diagram}
The Pontryagin product is given by the formula:
$$[P,a,f].[Q,b,g]=(-1)^{dim(P).|b|}[P\times Q,a\times b, \Upsilon_{\Omega M\times \Omega M}(f\times g)].$$
This product is associative but not commutative. The intersection
morphism is a morphism of algebras by commutativity of the diagram
above. $\square$
\\
\\
This morphism has been studied in details in \cite{FTV}, in
particular it is proved that the kernel of $I$ is nilpotent.

\subsection{Bordism and string topology} Let $\Omega^{SO}_*(X)$ be the bordism group of a
topological space $X$, we recall that it is isomorphic to the
bordism classes of singular oriented manifolds over $X$ (morphism
$f:M\rightarrow X$).
\\
We remark that $\Omega^{SO}_*(\mathcal LM)$ is also a $BV$-algebra
(all the constructions described above immediately adapt to
$\Omega^{SO}_*$).
\\
Let call a geometric class $[P,a,f]$ realizable if it is
equivalent to a class $[Q,1,g]$. This is equivalent to condition
of being in the image of the Steenrod-Thom morphism:
$$st:\Omega^{SO}_{n+d}(\mathcal LM)\rightarrow \mathbb H_n(M)$$
$$[M,f]\mapsto [M,1,f].$$
This morphism is clearly a morphism of $BV$-algebras and we have:

\subsubsection{\bf Proposition}{\it If $c\not \in Im(st)$ then
$I(c)=0$.}
\\
\\
{\bf Proof} This follows from the fact that if a geometric class
is not realizable (via the morphism $st$) it has the form
$[P,a,f]$ with $|a|>0$ and in this case $I([P,a,f])=0$. $\square$

\subsubsection{\bf Remark:} We recall that in general $st$ is neither injective nor surjective. However it is
surjective when $n+d<6$ (using Atiyah-Hirzebruch spectral sequence
one proves that it is an isomorphism for $n+d=0,1,2$) and it is
also surjective if we work over $\mathbb F_2$, over this field the
orientability condition in the definition of geometric homology is
unnecessary.

\subsection{String bracket}

\subsubsection{The string space} Let us consider the fibration:
$$S^1\rightarrow ES^1\rightarrow BS^1.$$
There exists a smooth model for this fibration, for $ES^1$ we take
$S^{\infty}$ the inductive limit of $S^n$. By \cite[chapter
X]{KrM} this is an Hilbert manifold modelled on $\mathbb
R^{(\mathbb N)}$. As $S^1$ acts freely and smoothly on
$S^{\infty}$ we have a $S^1$ fiber bundle of Hilbert manifolds:
$$S^1\rightarrow S^{\infty}\stackrel{\pi}{\rightarrow} \mathbb CP^{\infty}.$$
We get the $S^1$-fibration:
$$S^1\rightarrow \mathcal LM\times S^{\infty}\rightarrow \mathcal LM\times_{S^1}S^{\infty}.$$
The projection:
$$\mathcal LM\times S^{\infty}\rightarrow \mathcal LM$$
is a homotopy equivalence of Hilbert manifolds. As we know from
\cite{EE} that an homotopy equivalence between two separable
Hilbert manifolds is homotopic to a diffeomorphism, we deduce that
they are diffeomorphic.
\\
The space $\mathcal LM\times_{S^1}S^{\infty}$ is not an Hilbert
manifold because the action of $S^1$ on $\mathcal LM$ is not
smooth. Let call this space the string space of $M$.

\subsubsection{String homology} Let $\mathcal H_i$ be the homology group
$H'_{i+d}(\mathcal LM\times_{S^1}S^{\infty})$, this is the string
homology of $M$. In what follows we give explicit definitions of
the morphism $c$, $M$, $E$ of \cite[6]{CS}.
\\
\\
{\bf The morphism c.} Let $e\in H^2(\mathcal
LM\times_{S^1}S^{\infty})$ be the Euler class of the
$S^1$-fibration defined above:
$$c:\mathcal H_i\rightarrow \mathcal H_{i-2}$$
$$c([P,a,f])=[P,f^*(e)\cup a,f].$$
\\
\\
{\bf The morphism E.} This morphism is $E=\pi_*$:
$$E:\mathbb H_i(\mathcal LM)\rightarrow \mathcal H_i$$
$$E([P,a,f])=[P,a,\pi f].$$
\\
\\
{\bf The morphism M.} Let $[P,a,f]$ a geometric class in $\mathcal
H_i$. For a point $p\in P$, we can choose $(c_p,u)\in \mathcal
LM\times S^{\infty}$ that represents a class $[c_p,u]\in \mathcal
LM\times_{S^1}S^{\infty}$, as this choice is non-canonical we take
all the orbit of $(c_p,u)$ under the action of $S^1$. In this way
we produce a map
$$\tilde f: S^1\times P\rightarrow \mathcal
LM\times S^{\infty}.$$ Identifying $\mathcal LM\times S^{\infty}$
with $\mathcal LM$ one get a map:
$$M:\mathcal H_i\rightarrow \mathbb H_{i+1}(\mathcal LM)$$
$$M([P,a,f])=(-1)^{|a|}[S^1\times P,1\times a,\tilde f].$$
\\
\\
We have the following exact sequence, which is the Gysin exact
sequence associated to the $S^1$-fibration $\pi$:
$$\ldots\rightarrow \mathbb H_{i}(\mathcal LM)\stackrel{E}{\rightarrow}
\mathcal H_i \stackrel{c}{\rightarrow}\mathcal H_{i-2}
\stackrel{M}{\rightarrow}\mathbb H_{i-1}(\mathcal
LM)\rightarrow\ldots$$

\subsubsection{The bracket} The string bracket is given by the
formula:
$$[\alpha,\beta]=(-1)^{|a|}E(M(\alpha)\bullet M(\beta)).$$
Together with this bracket $(\mathcal H_*,[-,-])$ is a graded Lie
algebra of degree $(2-d)$ \cite[Th. 6.1]{CS}.

\subsection{Riemann surfaces operations} These operations are
defined by R. Cohen and V. Godin in \cite{CG} by means of Thom
spectra technology.
\\
Let $\Sigma$ be an oriented surface of genus $g$ with $p+q$
boundary components, $p$ incoming and $q$ outgoing. We fix a
parametrization of these components.
\\
Hence we have two maps:
$$\rho_{in}:\coprod_{p}S^1\rightarrow \Sigma,$$
and
$$\rho_{out}:\coprod_{q}S^1\rightarrow \Sigma.$$
If we consider the space of ${\bf H}^2$-maps ${\bf
H}^2(\Sigma,M)$, we get a Hilbert manifold and the diagram of
Hilbert manifolds:
$$\mathcal LM^{\times q}\stackrel{\rho_{out}}{\longleftarrow}
{\bf H}^2(\Sigma,M)\stackrel{\rho_{in}}{\longrightarrow}\mathcal
LM^{\times p}.$$ Let $\chi(\Sigma)$ be the Euler characteristic of
the surface. Using Sullivan's Chord diagrams it is proved in
\cite{CG} that the morphism
$${\bf H}^2(\Sigma,M)\stackrel{\rho_{in}}{\longrightarrow}\mathcal
LM^{\times p}$$ has  a homotopy model:
$${\bf H}^2(c,M)\stackrel{\rho_{in}}{\longrightarrow}\mathcal
LM^{\times p}$$ that is an embedding of Hilbert manifolds of
codimension $-\chi(\Sigma).d$. Hence by using the Gysin morphism
for Hilbert manifolds one can define the operation:
$$\mu_{\Sigma}:H'_{*}(\mathcal LM^{\times p})\stackrel{\rho_{in}^{!}}
{\longrightarrow}H'_{*+\chi(\Sigma).d}({\bf H}^2(\Sigma,M))
\stackrel{{\rho_{out}}_*}{\longrightarrow}H'_{*+\chi(\Sigma).d}(\mathcal
LM^{\times q}).$$ All these operations are parametrized by the
topological space of marked, metric chord diagrams $C\mathcal
F^{\mu}_{p,q}(g)$ \cite[sect1]{CG}. In the next section we
introduce a reduced version of Sullivan's chord diagrams and give
some algebraic properties of the associated operations.

\subsubsection{Sullivan's chord diagrams.}
In the preceding morphism $c\in\mathcal C\mathcal
F^{\mu}_{p,q}(g)$ is the {\it Sullivan's chord diagram} associated
to to the surface $\Sigma$. Let us recall the definition of
\cite{CG}:

\subsubsection{\bf Definition}{\it A metric marked Sullivan chord diagram $c$ of type (g;p,q) is
a metric fat graph:
\\
- a graph whose vertices are at least trivalent such that the
incoming edges are equipped with a cyclic ordering,
\\
- it has the structure of a compact metric space (details are
given in \cite[def. 1]{CG} and \cite[chapter 8]{I}),
\\
this fat graph represents a surface of genus $g$ with $p+q$
boundary components. The set of metric fat graphs is denoted by
$\mathcal Fat_{p,q}(g)$.
\\
The cycling ordering of the edges defines "boundary cycles". Pick
an edge and an orientation on it, then traverse it in the
direction of the orientation, this leads to a vertex, at this
vertex take the next edge coming from the cycling ordering and so
on. Then we get a cycle in the set of edges.
\\
The graph $c$ consists of a union of $p$ parametrized circles of
varying radii that represent the incoming boundary components,
joined at a finite number of points.
\\
Each boundary cycle has a marking. This marking correspond to the
starting point of a $S^1$-parametrization.}
\\
\\
Fat graphs (also called ribbon graphs) are a nice combinatoric
tool in order to study Riemann surfaces \cite{Penkava},
\cite{Penner} and \cite{strebel}.
\\
\\
Following a suggestion of A. Voronov rather than using the space
$C\mathcal F^{\mu}_{p,q}(g)$, we introduce the space of reduced
metric marked Sullivan chord diagrams denoted by
$\overline{C\mathcal F^{\mu}_{p,q}(g)}$. In a Sullivan diagram
there are ghost edges, these are the edges that lie on the $p$
disjoint circles. Hence there is a continuous map:
$$\pi:c\mapsto S(c)$$
that collapses each ghost edge to a vertex, let us notice that
$S(c)$ is also a metric fat graph but it is no more a Sullivan's
Chord diagram. And let $\overline{C\mathcal F^{\mu}_{p,q}(g)}$ be
the image of $\pi$. The space $\overline{C\mathcal
F^{\mu}_{p,q}(g)}$ has the following properties:

\subsubsection{\bf Proposition}{\it Let
$i:C\mathcal F^{\mu}_{p,q}(g)\rightarrow \mathcal Fat_{p,q}(g)$ be
the canonical inclusion then $i$ and $\pi$ are homotope.}

\subsubsection{\bf Proposition}{\it The space $\overline{C\mathcal F^{\mu}_{p,q}(g)}$ is a Prop.}
\\
\\
{\bf Proof.} The structure of Prop is obtained via a gluing
procedure by identifying an incoming boundary component with an
outgoing boundary component. This procedure is described at the
level of $C\mathcal F^{\mu}_{p,q}(g)$ in \cite{CG}. On unreduced
diagrams this procedure is not well defined because there is an
ambiguity coming from ghost edges, this is no more the case at the
level of $\overline{C\mathcal F^{\mu}_{p,q}(g)}$. $\square$

\subsubsection{\bf Proposition}{\it The suboperads $\overline{C\mathcal F^{\mu}_{*,1}(0)}$
and $\overline{C\mathcal F^{\mu}_{1,*}(0)}$ are isomorphic to the
$cacti$ operad.}
\\
\\
{\bf Proof.} This is immediate. Let us notice that any element in
the preimage $S^{-1}(\overline{c})$ where
$$S:C\mathcal F^{\mu}_{1,q}(0) \rightarrow \overline{C\mathcal F^{\mu}_{1,q}(0)}$$
is a chord diagram associated to the cactus $\overline{c}$ as
defined in \cite{K}. $\square$

\subsubsection{\bf Theorem}{\it For $q>0$ we have morphisms:
$$\mu_{n,p,q}(g):H'_n(\overline{C\mathcal F^{\mu}_{p,q}(g)})\rightarrow
Hom(H'_*(\mathcal LM^{\times p}),H'_{*+\chi(\Sigma).d+n}(\mathcal
LM^{\times q}))$$}
\\
\\
{\bf Proof} Let us give the construction of $\mu_{n,p,q}(g)$:
\\
consider an element of $H'_n(\overline{C\mathcal
F^{\mu}_{p,q}(g)})$ and suppose that is represented by a geometric
cycle $(S,\alpha,g)$ where $g:S\rightarrow \overline{\mathcal
C\mathcal F^{\mu}_{p,q}(g)}$. And let us define:
$$Map(g,M)=\{(s,f)/s\in S, f\in Map(g(s),M) \},$$
we also have maps:
$$\rho_{in}:Map(g,M)\rightarrow\mathcal LM^{\times p},$$
$$\rho_{out}:Map(g,M)\rightarrow\mathcal LM^{\times q}.$$
Using \cite[lemma3]{CG} we get an embedding of codimension $-\chi(\Sigma).d$:
$$\rho_{in}\times p:Map(g,M)\rightarrow \mathcal LM^{\times p}\times S$$
where $p$ is the canonical projection. If this map is an embedding
of Hilbert manifolds (we actually don't know if
$\overline{\mathcal C\mathcal F^{\mu}_{p,q}(g)}$ is a manifold),
we can form the following diagram:
\begin{diagram}
N_1\times\ldots\times N_p\times S&\lTo^{j}&N_1\ast\ldots\ast N_p\ast S\\
\dTo^{f_1\times\ldots \times f_p\times Id_S}&&\dTo^{\phi_g}\\
\mathcal LM^{\times p}\times S &\lTo_{\rho_{in}\times p}& Map(g,M)& \rTo_{\rho_{out}}& \mathcal LM^{\times q}.\\
\end{diagram}
We define $\mu_{n,p,q}$ in the following way:
\begin{diagram}
[N_1,\alpha_1,f_1]\otimes\ldots\otimes[N_p,\alpha_p,f_p]\\
\dTo\\
\pm[N_1\times\ldots\times N_p\times S,\alpha_1\times\ldots\times\alpha_p\times\alpha,f_1\times\ldots\times f_p\times id_S]\\
\dTo\\
\pm[N_1\ast\ldots\ast N_p\ast S,j^*(\alpha_1\times\ldots\times\alpha_p\times\alpha),\rho_{out}\phi_g].\\
\end{diagram}
Let us define the morphism in a less heuristic and more rigorous way. Rather than using an hypothetical embedding
of Hilbert manifolds let use the Thom collapse map of \cite[lemma5]{CG}:
$$\tau_g:\mathcal LM^{\times p}\times S\rightarrow (Map(g,M))^{\nu(g)}$$
where $\nu(g)$ is an open neighborhood of the embedding $\rho_{in}\times p$ and let $th_g$ be the Thom isomorphism:
$$th_g:H'_*((Map(g,M))^{\nu(g)})\rightarrow H'_{*+\chi(\Sigma).d}(Map(g,M)).$$
Then $\mu_{n,p,q}$ is defined by:
\begin{diagram}
[N_1,\alpha_1,f_1]\otimes\ldots\otimes[N_p,\alpha_p,f_p]\\
\dTo\\
\pm[N_1\times\ldots\times N_p\times S,\alpha_1\times\ldots\times\alpha_p\times\alpha,f_1\times\ldots\times f_p\times id_S]\\
\dTo\\
\pm[N_1\times\ldots\times N_p\times S,\alpha_1\times\ldots\times\alpha_p\times\alpha,
\tau_g(f_1\times\ldots\times f_p\times id_S)]\\
\dTo\\
\pm{\rho_{out}}_* th_g([N_1\times\ldots\times N_p\times S,\alpha_1\times\ldots\times\alpha_p\times\alpha,
\tau_g(f_1\times\ldots\times f_p\times id_S)]).\\
\end{diagram}

$\square$

\subsubsection{\bf Proposition}{\it $H'_*(\mathcal LM)$ is an algebra over the prop
$H'_*(\overline{C\mathcal F^{\mu}_{p,q}(g)})$}.
\\
\\
{\bf Proof} This follows immediately from \cite[Cor. 9]{CG}.
$\square$
\\
\\
This last result unifies some of the algebraic structure arising
in string topology, as corollaries we get:

\subsubsection{\bf Corollary}{\it If we fix $n=0$, then the action of homological degree $0$ string topology
operations induces on $H'_{*+d}(\mathcal LM)$ a structure of
Frobenius algebra without co-unit.}

\subsubsection{\bf Corollary}{\it When restricted to $H'_*(\overline{C\mathcal F^{\mu}_{p,1}(0)})$ on recovers
the BV-structure on $H'_{*+d}(\mathcal LM)$ induced by the string
product and the operator $\Delta$.}

\subsubsection{\bf A conjecture and its consequences.} We obtain homological string topology operations
which are 4-graded, 3 gradings being purely geometric ($p$, $q$ and $g$) and one grading purely homological.
\\
We do not know the homotopy type of $\overline{C\mathcal
F^{\mu}_{p,1}(0)}$. In fact there is a conjecture of R. Cohen and
V. Godin on the topology of the space $\mathcal CF^{\mu}_{p,q}$,
this conjecture states that this space is homotopy equivalent to
$B\mathcal M_g^{p+q}$, where $\mathcal M_g^{p+q}$ is the mapping
class group of a Riemann surface $S$ of genus $g$ with $p+q$
boundary components:
$$\mathcal M_q^{p+q}=\pi_0(Diff^+(S,\delta))$$
where the diffeomorphisms must preserve the boundary components
pointwise.
\\
Now suppose that this conjecture holds, and let fix $q=1$, since Tillmann's work \cite{Ti} we know that
the space $B\mathcal M_g^{p+1}$ is an operad $B\mathcal M$ and that this operad detects infinite loop spaces.
Let $\Gamma$ be the symmetric operad \cite{BE}, we have a map of operads:
$$\Gamma\rightarrow B\mathcal M.$$
This last fact
has a deep consequence in string topology because if we look at the singular homology of the free loop spaces
with coefficients in $\mathbb Z/p\mathbb Z$ this conjecture implies the existence of
"stringy Dyer-lashof" operations coming from the symmetries of the surfaces. From the operadic
nature of these operations they should satisfy Cartan and Adem relations, and these operations
should come from an $E_{\infty}$ structure \cite[chapter 1]{KM}.

\section{$n$-sphere topology.} Let fix an integer $n>1$, and
suppose that $M$ is a $n$-connected compact oriented smooth
manifold. Some results about the algebraic structure of the
homology of $n$-sphere spaces were announced in \cite[Th. 2.5]{V}.
We show how to recover a part of this structure.

\subsection{Sphere product.}Let us take $[P,a,f]\in
H'_{n_1+d}(\mathcal S_{n}M)$ and $[Q,b,g]\in H'_{n_2+d}(\mathcal
S_{n}M)$ two families of $n$-sphere. We can smoothen $f$ and $g$
and make them transverse to $\widetilde{\Delta}$, then we form the
pull-back $P\ast Q$.

\subsubsection{\bf Definition}{\it
Let $l:P\ast Q\rightarrow P$ and $r:P\ast Q\rightarrow Q$ be the
canonical maps, then we have the pairing:
$$-\bullet-:H'_{n_1+d}(\mathcal S_{n}M)\otimes H'_{n_2+d}(\mathcal S_{n}M)\longrightarrow
H'_{n_1+n_2+d}(\mathcal S_{n}M)$$
$$[P,a,f]\bullet[Q,b,g]=(-1)^{d.(|a|+|b|)+dim(P).|b|}[P\ast Q,l^*(a)\cup r^*(b),\Upsilon
\psi],$$ let call it the sphere product.}

\subsubsection{\bf Proposition}{\it The sphere product is associative
and commutative.}
\\
\\
{\bf Proof} The associativity and commutativity of the sphere
product follows from the associativity and commutativity of the
intersection product, the cup product and the fact that $\Upsilon$
is also associative and commutative up to homotopy. $\square$

\subsection{Constant spheres} We have a canonical embedding:
$$c:M\hookrightarrow \mathcal S_{n}M$$
$c$ induces a map:
$$c_*:H'_{*+d}(M)\rightarrow H'_{*+d}(\mathcal S_{n}M).$$
The morphism $c_*$ is clearly a morphism of commutative algebras.

\subsection{Intersection morphism}Let us recall that the map
$$ev_0:\mathcal S_{n}M\longrightarrow M$$
is a submersion (in fact this is a smooth fiber bundle of Hilbert
manifolds). Hence if we choose a base point $m\in M$ the fiber of
$ev_0$ in $m$ is the Hilbert manifold $\Omega^{n} M$ of
$n$-iterated based loops in $M$. Consider the morphism:
$$i:\Omega^n M\hookrightarrow \mathcal S_{n}M$$
this is an orientable morphism of codimension $d$.
\\
Let us describe the intersection morphism:
$$I=i^{!}:
H'_{*+d}(\mathcal S_{n}M)\rightarrow H'_*(\Omega M).$$ Let
$[P,a,f]\in H'_{n+d}(\mathcal S_{n}M)$ be a geometric class, one
can define $I([P,a,f])$ in the following way:
\\
\\
take the sphere product with $[c_m,1,c]$ where $c_m$ is a point
and $c:c_m\rightarrow \mathcal S_{n}M$ is the constant sphere at
the point $m$, then we set:
$$I([P,a,f])=(-1)^{d.|a|}[P\ast c_m,l^*(a),\psi].$$ As in the case $n=1$, we
remark we have:

\subsubsection{\bf Proposition}{\it The intersection morphism $I$ is a morphism of graded commutative
and associative algebras.}

\subsection{$3$-sphere topology} Let emphasis on the case $n=3$,
in this case we show the existence of an operator of degree $3$
acting on the homology. And we also obtain some results about the
$S^3$-equivariant homology of $3$-sphere spaces.

\subsubsection{The operator $\Delta_3$} The sphere $S^3$ acts on
$\mathcal S_3M$, let denote this action by $\Theta$. Hence, if we
consider a family of $3$-sphere in $M$
$$f:P\rightarrow \mathcal S_3M$$
we can build a new familly:
$$\Theta_f:S^3\times P\stackrel{id\times f}\longrightarrow
S^3\times \mathcal S_3M\stackrel{\Theta}{\rightarrow}\mathcal
S_{3}M.$$

\subsubsection{\bf Definition}{\it The operator $\Delta_3$ is given by the following formula:
$$\Delta_3:H'_{*}(\mathcal S_{3}M)\rightarrow H'_{*+3}(\mathcal S_{3}M)$$
$$[P,a,f]\mapsto (-1)^{|a|}[S^3\times P,1\times a,\Theta_f].$$}

\subsubsection{\bf Proposition}{\it The operator verifies:
$$\Delta_3^2=0.$$}
\\
\\
{\bf Proof} The proof is exactly the same as in the case of the
operator $\Delta$ and follows from the associativity of the cross
product and the nullity of the geometric class:
$$[S^3\times S^3,1\times 1, \mu].$$
$\square$

\subsubsection{$3$-sphere bracket} As in the case $n=1$, we have a smooth model of the
$S^3$-fibration:
$$S^3\rightarrow ES^3\rightarrow BS^3$$
which is given by the $S^3$ fiber bundle of Hilbert manifolds:
$$S^3\rightarrow S^{\infty}\rightarrow \mathbb HP^{\infty}.$$
We consider the
$S^3$-fibration:
$$S^3\rightarrow \mathcal S_3M\times ES^3\rightarrow\mathcal S_3M\times_{S^3} ES^3.$$
Let $\mathcal H_*^3=H_{*+d}(S_3M\times ES^3)$ and consider the
Gysin exact sequence of the fibration:
$$\ldots\rightarrow H'_{i+d}(\mathcal S_3M)\stackrel{E}{\rightarrow}
\mathcal H_i^3 \stackrel{c}{\rightarrow}\mathcal H_{i-4}^3
\stackrel{M}{\rightarrow}H'_{i+d-1}(\mathcal
S_3M)\rightarrow\ldots.$$

\subsubsection{\bf Definifion}{\it We define the $3$-sphere bracket by the formula:
$$[\alpha,\beta]=(-1)^{|\alpha|}E(M(\alpha)\bullet M(\beta)).$$}

\subsubsection{\bf Remark.} This bracket is anti-commutative.

\bigskip

\end{document}